\documentclass[11pt]{article}

\usepackage[cp1250]{inputenc}
\usepackage{multirow}
\usepackage{amsthm}
\usepackage{amsmath,amsfonts,amssymb,graphics,verbatim,
}

\newtheorem{theorem}{Theorem}

\newtheorem{lemma}[theorem]{Lemma}

\newtheorem{remark}[theorem]{Remark}

\def\N{{\mathbb N}}

\begin{document}
\title{A note on Wiener-Hopf factorization for
Markov Additive processes}

\author{Przemys\l aw Klusik
\footnote{University of Wroc\l aw, pl.\ Grunwaldzki 2/4, 50-384
Wroc\l aw, Poland, E-mail: przemyslaw.klusik@gmail.com}\quad
Zbigniew Palmowski\footnote{University of Wroc\l aw, pl.\
Grunwaldzki 2/4, 50-384 Wroc\l aw, Poland, E-mail:
zbigniew.palmowski@gmail.com}} \maketitle

\begin{abstract}
We prove the Wiener-Hopf factorization for
Markov Additive processes. We derive also Spitzer-Rogozin theorem
for this class of processes which serves for
obtaining Kendall's formula and Fristedt representation of
the cumulant matrix of the ladder epoch process.
Finally, we also obtain the so-called ballot theorem.
\end{abstract}

\section{Introduction}
The classical Wiener-Hopf factorization of a probability measure
was given by Spitzer (1964) and Feller (1970), and has a strong
connection to random walks. This result was generalized by Rogozin
(1966), Fristedt (1974) and other authors using approximation
based on discrete time skeletons. Greenwood and Pitman (1980) use
direct approach which relies on excursion theory for reflected
process. For details see Bertoin (1996) and Kyprianou (2006).
Presman (1969) and Arjas and Speed (1973) generalize Spitzer
identity into different direction, to the class of Markov Additive
processes in discrete time (see also Asmussen (2003) and Prabhu
(1998)). Later, Kaspi (1982) proves Wiener-Hopf factorization for
a continuous time parameter Markov Additive process, where
Markovian component has a finite state space and is ergodic. The
main weakness of the fluctuation identity given by Kaspi (1982) is
that they involve distribution of the inverse local time, which is
seldom explicitly known. Dieker and Mandjes (2006) investigate
discrete-time Markov Additive processes and use an embedding to
relate these to a continuous-time setting (see also Breuer (2008) and Rogers (1994)).

This paper presents Wiener-Hopf factorization for a special, but
none the less quite general, class of Markov Additive Processes
(MAP). For this class of processes we give short proof of
Wiener-Hopf factorization based on Markov property and additivity.
We also express the terms of Wiener-Hopf factorization directly in
terms of the basic data of the process. Finally,
we derive Spitzer-Rogozin theorem for this class of processes
which serves for obtaining  Kendall's formula and Fristedt
representation of the cumulant matrix of the ladder epoch process.
We also present the ballot theorem.

The paper is organized as follows. The Section \ref{prel}
introduces basic definitions, facts and properties related with MAPs.
In Section \ref{mainresults} we give the main results of this paper.
Finally in Section \ref{proofs} we prove all theorems.

\section{Preliminaries}\label{prel}
\subsection{Markov Additive Processes.}\label{scale.sec}

Before presenting main results we shall simply begin by defining
the class of processes we intend to work with and their
properties. Following Asmussen and Kella (2000) we consider a
process $X(t)$, where $X(t)=X^{(1)}(t)+X^{(2)}(t)$, and the
independent processes $X^{(1)}(t)$ and $X^{(2)}(t)$ are specified
by the characteristics: $q_{ij},G_{i},\sigma_i, a_i, \nu_i(dx)$
which we shall now define. Let $J(t)$ be a right-continuous,
ergodic, finite state space continuous time Markov chain, with
${\mathcal I}=\{1,\ldots,N\}$, and with the intensity matrix
${\mathbf Q}=(q_{ij})$. We denote the jumps of the process $J(t)$
by $\{T_i\}$ (with $T_0=0$). Let $\{U^{(i)}_n\}$ be
i.i.d. random variables with distribution function $G_{i}(\cdot)$.
Define the jump process by
$$
X^{(1)}(t)=\sum_{n\ge 1}\sum_{i\in\mathcal{I}} U^{(i)}_n
\mathbf{1}_{\{J(T_{n})=i,\ T_n\le t\}}.$$ For each $j\in \mathcal
I$, let $X^j(t)$ be a L\'{e}vy process with the L\'{e}vy-Khinchine
exponent:
\begin{eqnarray*}
&&-\log \mathbb{E}(\exp\{i\alpha X^j (1)\})
=\Psi_j(\alpha)\\&&\qquad\qquad
=-ia_j\alpha+\frac{\sigma_j^2\alpha^2}{2}+
\int_{-\infty}^\infty(1-e^{i\alpha y} +i\alpha |y| 1_{|y|\leq
1})\nu_j(dy),
\end{eqnarray*}
where $\int_{-\infty}^\infty (1\wedge |y|^2) \nu_j(dy)<\infty$. By
$X^{(2)}(t)$ we denote the process which behaves in law like
$X^i(t)$, when $J(t)=i$. We shall assume that the afore mentioned
class of MAPs is defined on a
probability space with probabilities $\left\{ \mathbb{P}_{i}:i\in \mathcal{%
I}\right\} $, where
$\mathbb{P}_i(\cdot)=\mathbb{P}(\cdot|J(0)=i)$, and
right-continuous natural filtration
$\mathbb{F=\{\mathcal{F}}%
_{t}:t\geq 0\mathbb{\}}$. In fact we can consider more general MAP
where additional jumps $U_n^{(i)}$ appearing during the change of
the state of $J(t)$ could also depend on the state $J(T_{n+1})$
(so called anticipative MAP). This could be done by considering
the vector state space $\mathcal{I}^2$ and the modified governing
Markov process $J$ on it. If each of the measures $\nu_i$ are supported
on $(-\infty,0)$ as well as the distributions of each $U^{(i)}$
then we say that $X$ is a {\it spectrally negative} MAP. These
definition and more concerning the basic characterization of MAPs
can be found in Chapter XI of Asmussen (2003).

\subsection{Time reversal}
Predominant in the forthcoming analysis will be the use of the bivariate
process $(\widehat{J},\widehat{X}),$ representing the process $(J,X)$ time
reversed from a fixed moment in the future when $J(0)$ has the stationary
distribution $\mathbf{\pi }$. For definitiveness, we mean
\begin{equation*}
\widehat{J}\left( s\right) =J( \left( t-s\right) ^{-} ) \text{ and }%
\widehat{X}\left( s\right) =X\left( t\right) -X\left( (t-s)^{-}\right) ,%
\text{ }0\leq s\leq t
\end{equation*}
under $\mathbb{P}_{\mathbf{\pi }}=\sum_{i\in \mathcal{I}}\pi _{i}\mathbb{P}%
_{i}.$ Note that $\widehat{X}$ is also
Markov Additive process. The characteristics of
$(\widehat{J},\widehat{X})$ will be indicated by using a hat over
the existing notation for the characteristics of $(J,X)$. Instead
of talking about the process $(\widehat{J},\widehat{X})$ we shall
also talk about the process $(J,X)$ under probabilities $\{\widehat{%
\mathbb{P}}_{i}:i\in \mathcal{I}\}$. Note also for future use,
following classical time reversed path analysis, for $y\geq 0$ and
$s\leq t$,
\begin{eqnarray}
\lefteqn{\mathbb{P}_{i}\left(\underline{G}(t)\in ds, -I(t)\in dy|J(t)=j\right)}\nonumber
\\&&=\widehat{\mathbb{P}}_{j}\left(\overline{G}(t)\in ds, S\left( t\right) -X(t)\in
dy|J(t)=i\right)\;,   \label{-I and S-X}
\end{eqnarray}
where $I(t)=\inf_{0\leq s\leq t}X(s)$, $S(t)=\sup_{0\leq s\leq t}X(s)$ and $\overline{G}(t)=\sup\{s<t: X(s)=S(s)\}$,
$\underline{G}(t)=\sup\{s<t: X(s)=I(s)\}.$ (A diagram may
help to explain the last identity).

\subsection{Ladder height process}
We start from recalling the representation of the
local time given in Kaspi (1982) in formula (3.21). For MAP we say that state
$i\in\{1,\ldots,N\}$ is regular when
$\mathbb{P}_{i}(R=0)=1$
for
$R=\inf\{t\ge 0: t\in \overline{\mathcal{M}}\},$ where
$\overline{\mathcal{M}}$ is a closure of $\mathcal{M}=\{t\ge 0:
X(t)=S(t)\}$. Denote by $\{U_n\}$ the stopping times at which
$R(t-)=0$ and $R(t)>0$ for the $R(t) =\inf\{s> t:
S(t)=X(t)\}-t$ and $J(t)$ is irregular. Denote
$$S_n=\left\{\begin{array}{ll}
U_n & \mbox{on $X(U_n)=S(U_n)$,}\\
\infty& \mbox{otherwise.}
\end{array}\right.
$$
By Theorem 3.28 of Kaspi (1982) (see also Maissoneuve (1975)), for
the MAP we can define the ladder height process:
$$\{(L^{-1}(t),
H(t)=X(L^{-1}(t)), J(L^{-1}(t))), t\geq 0\}$$ choosing the local
time:
\begin{equation}\label{localtime}L(t)=L^{\rm c}(t) +\sum_{S_n<t} \lambda
(J(U_n)) \mathbf{e}_1^{(n)},\end{equation} where $L^{\rm c}(t)$ is
a continuous additive process that increases only on $\mathcal{M}$
and $\mathbf{e}_1^{(n)}$ are independent exponential random
variables with intensity $1$,
$$\lambda(i)=\mathbb{E}_i\left[(1-e^{-R}\right].
$$
Obviously to make functional (\ref{localtime}) measurable we
enlarge probability space to include the exponential random
variables. One can easily verify that $(L^{-1}(t), H(t),
J(L^{-1}(t)))$ is again (bivariate) MAP (see Kaspi (1982, p.
185)). For each moment of time we can define the excursion:
$$\epsilon_t(s)=
\left\{\begin{array}{ll} X(L^{-1}(t-)+s) - X(L^{-1}(t-))&
\mbox{for $L^{-1}(t-)<L^{-1}(t)$}\\
\partial& \mbox{otherwise,}
\end{array}
\right. $$ where $\partial$ is a cementary state. Let
$\zeta(\epsilon_t)=L^{-1}(t)-L^{-1}(t-)$ be the length of the
excursion. From (\ref{localtime}) it follows that the excursion
process $\{(t,\epsilon_t),\;t\geq 0\}$ is (possibly stopped at the
first excursion with infinite length) marked Cox point process
with the intensity $n(J(L^{-1}(t-)), d\epsilon)$ depending on the
state process $J(L^{-1}(t-))$. Denote by $\mathcal{E}$ the
$\sigma$-field on the excursion state space.

\subsection{Spectrally negative Markov Additive process}
Letting
${\mathbf{Q}}\circ{\widetilde{\mathbf{G}}}(\alpha)=(q_{ij}\widetilde{G}_{i}(\alpha)
)$, where $\widetilde{G}_{i}(\alpha)= \mathbb{E}\left(\exp(\alpha
U^{(i)})\right)$, for {\it spectrally negative} MAP we can define
{\it cumulant generating matrix} (cgm) of MAP $X(t)$:
\begin{equation}{\mathbf{F}}(\alpha)={\mathbf{Q}}\circ{\widetilde{\mathbf{G}}}(\alpha)
+{\rm diag}(\psi_1(\alpha),\ldots,\psi_N(\alpha)),\qquad \alpha
\in \mathbb{R}_+\ . \label{defn:F}
\end{equation}
Perron-Frobenius theory identifies $\mathbf{F}%
\left( \alpha \right) $ as having a real-valued eigenvalue with
maximal absolute value which we shall label $\kappa \left( \alpha
\right) .$ The corresponding left and right $1\times N$
eigenvectors we label $\mathbf{v}\left( \alpha \right) $ and
$\mathbf{h%
}\left( \alpha \right) $, respectively. In this text we shall
always write vectors in their horizontal form and use the usual
$^{\text{T}}$ to mean transpose. Since $\mathbf{v}\left( \alpha
\right) $ and $\mathbf{h}\left( \alpha \right) $ are given up to
multiplying constants, we are free to normalize them such that
\begin{equation*}
\mathbf{v}\left( \alpha \right) \mathbf{h}\left( \alpha \right) ^{\text{T}}=1%
\text{ and }\mathbf{\pi h}\left( \alpha \right) ^{\text{T}}=1\;.
\end{equation*}
 Note also that $\mathbf{h}\left( 0\right) =\mathbf{e}$, the $1\times N$
vector consisting of a row of ones. We shall write $h_{i}\left(
\alpha \right) $ for the $i$-th element of $\mathbf{h}\left(
\alpha \right) .$ The eigenvalue $\kappa \left( \alpha \right) $
is a convex function (this can also be easily verified) such that
$\kappa \left( 0\right) =0$ and $\kappa ^{\prime }\left( 0\right)
$ is the asymptotic drift of $X$ in the sense that for each $i\in
\mathcal{I}$ we have $\lim_{t\uparrow \infty }$ $\mathbb{E}
(X(t)|J(0)=i, X(0)=x)/t=\kappa ^{\prime }\left( 0\right) $. For
the right inverse of $\kappa$ we shall write $\Phi$.

It can be checked that under the following Girsanov change of
measure
\begin{equation}
\left. \frac{d\mathbb{P}_{i}^{\gamma }}{d\mathbb{P}_{i}}\right| _{%
\mathcal{F}_{t}}:=e^{\gamma X\left( t\right)-\kappa \left( \gamma
\right) t}\frac{h_{J\left( t\right) }\left( \gamma \right)
}{h_{i}\left( \gamma \right) },\text{ for }\gamma \text{ such that
}\kappa \left( \gamma \right) <\infty   \label{alaGirsanov}
\end{equation}
the process $(X,\mathbb{P}_{i}^{\gamma })$ is again a spectrally
negative MAP whose intensity matrix $\mathbf{F}_{\gamma }\left(
\alpha \right) $ is well defined and finite for $\alpha \geq
-\gamma $. Generally for all quantities calculated for
$\mathbf{P}^{\gamma }$ we will add subscript $\gamma$. Further, if
$\mathbf{F}_{\gamma }\left( \alpha \right) $ has largest
eigenvalue $\kappa _{\gamma }\left( \alpha \right) $ and
associated right eigenvector $\mathbf{h}_{\gamma }\left( \alpha
\right) $, the triple $\left( \mathbf{F}_{\gamma }\left( \alpha
\right) ,\kappa _{\gamma }\left( \alpha \right)
,\mathbf{h}_{\gamma }\left( \alpha \right) \right) $ is related to
the original triple $\left(
\mathbf{F}\left( \alpha \right) ,\kappa \left( \alpha \right) ,\mathbf{h}%
\left( \alpha \right) \right) $ via
\begin{equation}
\mathbf{F}_{\gamma }\left( \alpha \right) =\mathbf{\Delta}
_{\mathbf{h}}\left( \gamma
\right) ^{-1}\mathbf{F}\left( \alpha +\gamma \right) \mathbf{\Delta} _{\mathbf{h}%
}\left( \gamma \right) -\kappa \left( \gamma \right) \mathbf{I}\text{ and }%
\kappa _{\gamma }\left( \alpha \right) =\kappa \left( \alpha
+\gamma \right) -\kappa \left( \gamma \right)\;,
\label{meas-change-effect}
\end{equation}
where $\mathbf{I}$ is the $N\times N$ identity matrix and
\begin{equation*}
\mathbf{\Delta} _{\mathbf{h}}\left( \gamma \right)
:=\text{diag}\left( h_{1}\left( \gamma \right) ,...,h_{N}\left(
\gamma \right) \right) .
\end{equation*}

Similarly, the time reversed process $\widehat{X}(t)$ is the
spectrally negative MAP with the characteristics
$\widehat{\mathbf{F}},$ $\widehat{\mathbf{h}}$, $\widehat{\kappa
}$. To relate them to the original ones, recall that the intensity
matrix of $\widehat{J}$ must satisfy
\begin{equation*}
\widehat{\mathbf{Q}}=\mathbf{\Delta} _{\mathbf{\pi
}}^{-1}\mathbf{Q}^{\text{T}}\mathbf{\Delta} _{\mathbf{\pi }}\;,
\end{equation*}
where $\mathbf{\Delta} _{\mathbf{\pi }}$ is the diagonal matrix
whose entries are given by the vector $\mathbf{\pi }$. Hence
according to (\ref{defn:F}) we find that:
\begin{equation*}
\widehat{\mathbf{F}}\left( \alpha \right) =\mathbf{\Delta} _{\mathbf{\pi }}^{-1}%
\mathbf{F}\left( \alpha \right) ^{\text{T}}\mathbf{\Delta}
_{\mathbf{\pi }}.
\end{equation*}
Moreover, $\widehat{\kappa }\left( \alpha \right) = \kappa \left(
\alpha \right) $ and
$\mathbf{\Delta} _{\mathbf{\pi }}\widehat{\mathbf{h}}\left( \alpha \right) ^{\text{T}%
}=\mathbf{v}\left( \alpha \right) ^{\text{T}}$ (see Kyprianou and
Palmowski (2008) for details).

The spectrally negative MAP is easier to analyze since its ladder
height process $(L^{-1}(t), H(t), J(L^{-1}(t)))$ has explicit
matrix cumulant generating matrix $\mathbf{\Xi}(q,\alpha)$. Let
\[
\tau^+_a:= \inf\{t\geq 0: X(t)\geq a\},
\]
where $a\geq 0$. Denote the generator of the Markov process
$\{J(\tau^+_a),a\geq 0\}$ by $\mathbf{\Lambda} (q)$ on
$\mathbb{P}^{\Phi(q)}$. It solves equation:
\begin{equation}\label{Lambdaeq}
\mathbf{F}_{\Phi(q)}(-\mathbf{\Lambda}(q))=\mathbf{0};\end{equation}
see Pistorius (2005) and Ivanovs et al. (2008). Note that ladder
height process can be identified as $\{(\tau_a^+,
X(\tau^+_a)=a,J(\tau_a^+)),\;a\ge 0\}$. It is a bivariate Markov
additive process with the cumulant generating matrix:
\begin{eqnarray}
\mathbf{\Xi}(q,\alpha)& =&\mathbf{\Delta}_{\mathbf{h}}(\Phi(q))
(\Phi(q)\mathbf{I}-\mathbf{\Lambda}(q))
\mathbf{\Delta}_{\mathbf{h}}(\Phi(q))^{-1}+\alpha\nonumber\\&=&
\mathbf{\Delta}_{\mathbf{h}}(\Phi(q))
((\Phi(q)+\alpha)\mathbf{I}-\mathbf{\Lambda}(q))
\mathbf{\Delta}_{\mathbf{h}}(\Phi(q))^{-1}\label{Xi}
\end{eqnarray}
for $\alpha,\;q>0$.
Above could be
deduced from the equalities:
\begin{equation}\label{Xipred}e^{-\mathbf{\Xi}(q,\alpha)a}=\mathbf{E}\left(e^{-q\tau_a^+-\alpha X(\tau_a^+)};
J(\tau^+_a)\right)=\mathbf{E}\left(e^{-q\tau_a^+-\alpha a};
J(\tau^+_a)\right)\end{equation}
and the Theorem 1 of Kyprianou and Palmowski (2008)
stating that
\begin{eqnarray} \lefteqn{\mathbf{E}(e^{-\xi \tau_x^+
};\tau_x^+<\mathbf{e}_q; J(\tau^+_x))= \mathbf{E}\left(
e^{-(q+\xi)\tau _{x}^{+}}1_{\mathbf{(}\tau _{x}^{+}<\infty
\mathbf{)}} ; J(\tau^+_x)\right)\nonumber}\\&& =
\mathbf{\Delta}_{\mathbf{h}}(\Phi(q+\xi))
e^{-(\Phi(q+\xi)\mathbf{I}-\mathbf{\Lambda}(q+\xi))x}
\mathbf{\Delta}_{\mathbf{h}}(\Phi(q+\xi))^{-1}\nonumber\\
&& = \exp\left\{-\mathbf{\Delta}_{\mathbf{h}}(\Phi(q+\xi))
(\Phi(q+\xi)\mathbf{I}-\mathbf{\Lambda}(q+\xi))
\mathbf{\Delta}_{\mathbf{h}}(\Phi(q+\xi))^{-1}x\right\}.
\label{up-crossing}
\end{eqnarray}

\section{Main results}\label{mainresults}
As much as possible, from now on, we shall prefer to work with
matrix notation. For a random variable $Y$ and (random) time
$\tau$, we shall understand $\mathbf{E}(Y; J(\tau))$ to be the
matrix with $(i,j)$-th elements $\mathbb{E}_{i}(Y; J(\tau)=j)$.
For an event, $A$, $\mathbf{P}(A; J(\tau))$ will be understood in
a similar sense. Here and throughout we work with the definition
that $\mathbf{e}_q$ is random variable which is exponentially
distributed with mean $1/q$ and independent of $(J,X)$. Let
$\mathbf{I}_{ij}(q) = \mathbb{P}_{i,0}(J(\mathbf{e}_q)=j)$, in
other words $$\mathbf{I}(q)=q(q\mathbf{I}-\mathbf{Q})^{-1}.$$ From
now on we assume that {\it none of the processes $X^i$ are
downward subordinators and compound Poisson process}. To include
compound Poisson process $X^{(i)}(t)$ in Theorem \ref{WH.main}(i)
on the event $\{J(\underline{G}(\mathbf{e}_q))=i\}$ it is
necessary to work with the new definition
$\underline{G}(t)=\inf\{s<t: X(s)=I(s)\}$ instead the previous
one.

\begin{theorem}\label{WH.main}
(i)  For a general MAP the random vectors
$(S(\mathbf{e}_q),\overline{G}(\mathbf{e}_q))$ and
$(S(\mathbf{e}_q)-X(\mathbf{e}_q),\mathbf{e}_q-\overline{G}(\mathbf{e}_q))$
are conditionally on $J(\overline{G}(\mathbf{e}_q))$ independent.
Hence for $\alpha \in \mathbf{R}$, $\xi\geq 0$,
\begin{eqnarray}
&&
\mathbf{E}\left[ e^{i\alpha X(\mathbf{e}_q)-\xi \mathbf{e}_q};
J(\mathbf{e}_q)\right] \label{th3i1b}\\&&= \mathbf{E}\left[
e^{i\alpha S(\mathbf{e}_q)-\xi \overline{G}(\mathbf{e}_q)};
J(\overline{G}(\mathbf{e}_q))\right]
\mathbf{\Delta}_{\mathbf{\pi}}^{-1}\widehat{\mathbf{E}}\left[
e^{i\alpha I(\mathbf{e}_q)-\xi \underline{G}(\mathbf{e}_q)};
J(\underline{G}(\mathbf{e}_q))\right]^{\text T}
\mathbf{\Delta}_{\mathbf{\pi}}\nonumber
\end{eqnarray}
and
\begin{eqnarray}
&& \mathbf{E}\left[ e^{i\alpha X(\mathbf{e}_q)-\xi \mathbf{e}_q};
J(\mathbf{e}_q)\right]\label{th3i1}\\&&= \mathbf{E}\left[
e^{i\alpha I(\mathbf{e}_q)-i\xi \underline{G}(\mathbf{e}_q)};
J(\underline{G}(\mathbf{e}_q))\right]
\mathbf{\Delta}_{\mathbf{\pi}}^{-1}\widehat{\mathbf{E}}\left[
e^{i\alpha S(\mathbf{e}_q)-i\xi \overline{G}(\mathbf{e}_q)};
J(\overline{G}(\mathbf{e}_q))\right]^{\text T}
\mathbf{\Delta}_{\mathbf{\pi}}\;.\nonumber
\end{eqnarray}

(ii)\ \ For the spectrally negative MAP and $\alpha,\xi\geq 0$,
\begin{equation}
\mathbf{E}\left[ e^{-\alpha S(\mathbf{e}_q)-\xi
\overline{G}(\mathbf{e}_q)}; J(\overline{G}(\mathbf{e}_q)) \right]
=\mathbf{\Xi}(q+\xi,\alpha)^{-1}{\rm diag}
\left(\mathbf{\Xi}(q,0)\mathbf{I}(q)\mathbf{e}^{\text
T}\right),\label{th3ii1b}
\end{equation}
\begin{equation}
\mathbf{E}\left[ e^{-\alpha S(\mathbf{e}_q)-\xi
\overline{G}(\mathbf{e}_q)}; J(\mathbf{e}_q) \right]=
\mathbf{\Xi}(q+\xi,\alpha)^{-1}\mathbf{\Xi}(q,0)\mathbf{I}(q),\label{th3ii1}
\end{equation}
\begin{eqnarray} \lefteqn{\mathbf{E}\left[
e^{\alpha I(\mathbf{e}_q)-\xi \underline{G}(\mathbf{e}_q)};
J(\underline{G}(\mathbf{e}_q)) \right]}\label{th3ii2b}\\&&=
q\left((q+\xi)\mathbf{I}-\mathbf{F}(\alpha)\right)^{-1}
\mathbf{\Delta}_{\mathbf{\pi}}^{-1}\widehat{\mathbf{\Xi}}(q+\xi,-\alpha)^{\text
T}{\rm diag}
\left(\widehat{\mathbf{\Xi}}(q,0)\widehat{\mathbf{I}}(q)\mathbf{e}^{\text
T}\right)^{-1}\mathbf{\Delta}_{\mathbf{\pi}},\nonumber\end{eqnarray}
\begin{eqnarray} \lefteqn{\mathbf{E}\left[
e^{\alpha I(\mathbf{e}_q)-\xi \underline{G}(\mathbf{e}_q)};
J(\mathbf{e}_q) \right]}\nonumber\\&&=
q\left((q+\xi)\mathbf{I}-\mathbf{F}(\alpha)\right)^{-1}
\mathbf{\Delta}_{\mathbf{\pi}}^{-1}\widehat{\mathbf{\Xi}}(q+\xi,-\alpha)^{\text
T}\left[\widehat{\mathbf{\Xi}}(q,0)^{-1}\right]^{\text
T}\mathbf{\Delta}_{\mathbf{\pi}} \;.\label{th3ii2}\end{eqnarray}
\end{theorem}

\begin{remark}\label{remarkreversedWH}
\rm Applying Theorem \ref{WH.main}(i) to the reversed process
derives similar conclusion for the infimum functional.
Namely, processes $\{(X(t),J(t)),0\le t<
\underline{G}(\mathbf{e}_q)\}$ and
$\{(X(\underline{G}(\mathbf{e}_q)+t)-X(\underline{G}(\mathbf{e}_q)),
J(\underline{G}(\mathbf{e}_q)+t)),t\ge 0\}$ are conditionally on
$J(\underline{G}(\mathbf{e}_q))$ independent.
\end{remark}

\begin{remark} \rm
For $N=1$ (hence $\mathbf{\Lambda}(q)=0$, $\mathbf{I}(q)=1$) above
theorem gives well-known identities for the spectrally negative
L\'{e}vy process:
\begin{eqnarray*}\mathbb{E}\left[e^{-\alpha S(\mathbf{e}_q)-\xi
\overline{G}(\mathbf{e}_q)}\right]&=&
\frac{\Phi(q)}{\Phi(q+\xi)+\alpha},\\ \mathbb{E}\left[e^{\alpha
I(\mathbf{e}_q)-\xi \underline{G}(\mathbf{e}_q)}\right] &=&
\frac{q(\Phi(q+\xi)-\alpha)}{\Phi(q)(q+\xi-\psi(\alpha))},
\end{eqnarray*}
where $\psi(\theta)=-\Psi(-i\theta)$ is a Laplace exponent of $X$.
Finally, for $\xi=0$ above theorem gives already known identity
for spectrally negative MAP (see Kyprianou and Palmowski (2008)):
\begin{eqnarray}
\lefteqn{\mathbf{E}\left( e^{\alpha I(\mathbf{e}_{q})} ;
J(\mathbf{e}_q)\right) ^{\text{T}} \left( \mathbf{F}\left( \alpha
\right) -q\mathbf{I}\right)
^{\text{T}} } \nonumber \\
&=&q\mathbf{\Delta }_{\mathbf{v}}\mathbf{(}\Phi \left( q\right) )
[ \alpha
(\Phi(q)\mathbf{I}    -\widehat{\mathbf{\Lambda}}(q)  )^{-1} -
\mathbf{I}] \mathbf{\Delta} _{\mathbf{v}}\mathbf{(}\Phi \left(
q\right) )^{-1}\;, \label{key-identity}
\end{eqnarray}
which was derived using Asmussen-Kella martingale.
\end{remark}

We prove also the
following counterpart of Spitzer-Rogozin version of Wiener-Hopf
factorization and the Fristedt theorem: 

\begin{theorem}\label{rogozin}
Assume that the matrix $\mathbf{E}\exp\{i\alpha X(1)\}$ has
distinct eigenvalues and that for any $t,s\geq 0$:
\begin{eqnarray}
\lefteqn{\mathbf{E}\left[
e^{i\alpha X(t)}\mathbf{1}_{\{X(t)\geq 0\}};
J(t) \right]\;\mathbf{E}\left[
e^{i\alpha X(s)}\mathbf{1}_{\{X(s)< 0\}};
J(s) \right]}\nonumber\\&&=\mathbf{E}\left[
e^{i\alpha X(s)}\mathbf{1}_{\{X(s)< 0\}};
J(s) \right]\;\mathbf{E}\left[
e^{i\alpha X(t)}\mathbf{1}_{\{X(t)\geq 0\}};
J(t) \right].\label{commute}\end{eqnarray}
Then
\begin{eqnarray*} \lefteqn{\mathbf{E}\left[
e^{-\alpha S(\mathbf{e}_q)-\xi \overline{G}(\mathbf{e}_q)};
J(\mathbf{e}_q) \right]}\\&&=\exp\left\{\int_0^\infty
dt\int_{[0,\infty)}\left( e^{-\xi t-\alpha x}-1\right)t^{-1}e^{-q
t}\mathbf{P}(X(t)\in dx; J(t))\right\}\mathbf{I}(q)
\end{eqnarray*}
and
\begin{eqnarray*} \lefteqn{\mathbf{E}\left[
e^{\alpha I(\mathbf{e}_q)-\xi\underline{G}(\mathbf{e}_q)};
J(\mathbf{e}_q) \right]}\\&&=\exp\left\{\int_0^\infty
dt\int_{(-\infty,0)}\left( e^{-\xi t+\alpha x}-1\right)t^{-1}e^{-q
t}\mathbf{P}(X(t)\in dx; J(t) )\right\}\mathbf{I}(q)\;.
\end{eqnarray*}
\end{theorem}

The assumption (\ref{commute}) is satisfied for example for
Markov modulated Brownian motion $X(t)=\sigma(J(t))B(t)$, where $\sigma$ is a
positive function.

The following generalization of the Kendall's identity and the ballot theorem also hold.
\begin{theorem}\label{kendall}
Consider the spectrally negative Markov Additive process $X(t)$. If
there exist distinct $q_1, q_2,\ldots,q_N$ such that vectors
$\mathbf{h}(\Phi(q_1)), \mathbf{h}(\Phi(q_2)),$ $\ldots,
\mathbf{h}(\Phi(q_N))$ are independent, then
$$t\mathbf{P}(\tau_x^+\in dt; J(t))dx=x\mathbf{P}(X(t)\in dx;J(t))dt\;.$$
\end{theorem}
\begin{theorem}\label{ballot}
Let $$X(t)=ct-\sigma(t),$$ where $\{\sigma(t),t\ge 0\}$ is a
Markov Addidive Subordinator without drift component. Under the assumptions
of the Theorem \ref{kendall}, the following identity holds:
$$\mathbf{P}(X(t)\in dx, I(t)=0;J(t))=\frac{x}{ct}\mathbf{P}(X(t)\in dx;J(t))\;.$$
\end{theorem}
One can straightforward check that the assumptions
of the Theorem \ref{kendall} are satisfied e.g. for $X(t)=ct-J(t)N(t)$, where
$N(t)$ is a Poisson process and $J(t)$ is a two-state birth-death process.

In total theorems given here might be seen as a fundamental of the
fluctuation theory for the (spectrally negative) MAP and might serve
for the deriving counterparts of well-known identities for the
L\'{e}vy processes.

\section{Proofs}\label{proofs}
\subsection{Proof of Theorem \ref{WH.main}}\label{ii}
(i)  Sampling MAP process $(X(t), J(t))$ up to exponential random time
$\mathbf{e}_q$ corresponds to the sampling the marked Cox point process (double Poisson point process)
of the excursions up to time
$L(\mathbf{e}_q)$.
Moreover, since
conditioning on realization of the process $J(t)$ the point
process $(t,\epsilon_t)$ is a non-homogeneous marked Poisson process, we
know that, conditioning on $J(L^{-1}(\sigma^A-))$
for
$$\sigma^A=\inf\{t\geq 0: \epsilon_t \in A\},$$
the point process
$\{(t,\epsilon_t), t<\sigma^A\}$ is independent of
$\epsilon_{\sigma^A}$.

Consider now
$$\sigma_1=\inf\left\{t\geq 0:
\int_0^{L^{-1}(t)}\mathbf{1}_{\{X(s)=S(s)\}}\;ds>\mathbf{e}_q\right\}$$
and
$$\sigma_2=\inf\{t\geq 0: \zeta(\epsilon_t)>\mathbf{e}_q^t\},$$
where $\mathbf{e}_q^t$ is the independent exponential random
variable with intensity $q$ if $\epsilon_t\neq\partial$ and
$\mathbf{e}_q^t=\partial$ otherwise. Note that $\sigma_2$
is $\sigma^A$ for the $A=\{\zeta(\epsilon)>\mathbf{e}_q\}$. If
$\sigma_2<\sigma_1$, then conditioning on
$J(L^{-1}((\sigma_1\wedge\sigma_2)-))=J(L^{-1}(\sigma_2-))$ the
process
\begin{equation}\label{ecursionpast}
\{(t,\epsilon_t),\quad
t<\sigma_1\wedge\sigma_2\quad\mbox{and}\quad \epsilon_t \neq
\partial \}
\end{equation}
is independent of
$\epsilon_{\sigma_2}=\epsilon_{\sigma_1\wedge\sigma_2}$. If
$\sigma_1<\sigma_2$, then
$\epsilon_{\sigma_1}=\epsilon_{\sigma_1\wedge\sigma_2}=\partial$
and is also independent of the process (\ref{ecursionpast}). Hence
conditioning on $J(L^{-1}(\sigma_1\wedge\sigma_2-))$ the excursion
$\epsilon_{\sigma_1\wedge\sigma_2}$ is independent of the process
(\ref{ecursionpast}).
Note also that
\begin{equation}\label{represen}\overline{G}(\mathbf{e}_q)=L^{-1}((\sigma_1\wedge\sigma_2)-),\qquad
S(\mathbf{e}_q)=H((\sigma_1\wedge\sigma_2)-)\end{equation} and the
last excursion $\epsilon_{\sigma_1\wedge\sigma_2}$ occupies the
final $\mathbf{e}_q-\overline{G}(\mathbf{e}_q)$ units of time in
the interval $[0,\mathbf{e}_q]$ and reaches the depth
$X(\mathbf{e}_q)-S(\mathbf{e}_q)$. This completes the proof of the first part of the
Theorem \ref{WH.main}(i). Note that
$(\mathbf{e}_q-\overline{G}(\mathbf{e}_q, X(\mathbf{e}_q)-S(\mathbf{e}_q))$ has the same law like
$(\widehat{\underline{G}}(\mathbf{e}_q), \widehat{I}(\mathbf{e}_q))$. The second part of of
the Theorem \ref{WH.main}(i) follows now from the first part
applied for the reversed process.

(ii) To prove Theorem \ref{WH.main}(i) we follow Bertoin (2000).
Fix $n\in\N$ setting
\begin{equation}\label{in}
i_n^+=[nS(\mathbf{e}_q)]/n,\end{equation} where $[\cdot]$ stands
for integer parts. Applying the strong Markov property at time
$\tau_{k/n}^+$ and using (\ref{up-crossing}) yields:
\begin{eqnarray}
\lefteqn{\mathbf{E}\left[ e^{-\xi\tau_{i_n^+}^+-\alpha
i_n^+};J(\mathbf{e}_q)\right]}\nonumber\\
&&= \sum_{k=0}^\infty \mathbf{E}\left[e^{-\xi\tau_{k/n}^+-\alpha
k/n };k/n\leq S(\mathbf{e}_q)<(k+1)/n;J(\mathbf{e}_q)\right]\nonumber\\
&&= \sum_{k=0}^\infty
e^{-\mathbf{\Xi}(q+\xi,\alpha)k/n}\mathbf{P}\left(\mathbf{e}_q<\tau_{1/n}^+;
J(\mathbf{e}_q) \right)\label{basicidentitywha}
\\&& =
\left[n(\mathbf{I}-e^{-\mathbf{\Xi}(q+\xi,\alpha)1/n})
\right]^{-1}n\mathbf{P}\left(\mathbf{e}_q<\tau_{1/n}^+;
J(\mathbf{e}_q) \right).\label{basicidentitywh}
\end{eqnarray}
Taking $n\to\infty$ we have $i_n^+\to S(\mathbf{e}_q)$ and
$\tau_{i_n^+}^+\to \overline{G}(\mathbf{e}_q)$ and hence the left
hand side of above equation converges by the dominated convergence
theorem. Thus also right hand side converges. Note that for any
matrix $\mathbf{A}$
\begin{equation}\label{asaa}\mathbf{I}-e^{-\mathbf{A}/n}=\frac{1}{n}
\mathbf{A}+{\text o}(1/n).\end{equation} Thus
$$
\mathbf{E}\left[e^{-\xi\overline{G}(\mathbf{e}_q)-\alpha
S(\mathbf{e}_q) };J(\mathbf{e}_q)\right]
=\mathbf{\Xi}(q+\xi,\alpha)^{-1}\mathbf{B}
$$
for some matrix $\mathbf{B}$ using fact that matrix
$(\Phi(q+\xi)+\alpha)\mathbf{I}-\mathbf{\Lambda}(q+\xi)$ is
invertible for $q>0$. Taking $\xi=\alpha=0$ we obtain
$$\mathbf{B}=\mathbf{\Xi}(q,0)
\mathbf{I}(q)
$$
which completes the proof of (\ref{th3ii1}). Similarly, from
(\ref{basicidentitywha}),
\begin{eqnarray*}
\lefteqn{\mathbb{E}_i\left[ e^{-\xi\tau_{i_n^+}^+-\alpha
i_n^+};J(\tau_{i_n^+}^+)=j\right]}\nonumber\\
&&= \sum_{k=0}^\infty \mathbb{E}_i\left[e^{-\xi\tau_{k/n}^+-\alpha
k/n };\tau^+_{k/n}<\mathbf{e}_q;J(\tau_{k/n}^+)=j\right]\mathbb{P}_j\left(\mathbf{e}_q<\tau_{1/n}^+\right)\nonumber\\
&&=\left[\mathbf{E}\left[ e^{-\xi\tau_{i_n^+}^+-\alpha
i_n^+};J(\mathbf{e}_q)\right]
\mathbf{P}\left(\mathbf{e}_q<\tau_{1/n}^+; J(\mathbf{e}_q)
\right)^{-1}\right]_{ij}\mathbb{P}_j\left(\mathbf{e}_q<\tau_{1/n}^+\right)
\end{eqnarray*}
and therefore \begin{equation}\label{pathdecomp2}
\mathbf{E}\left[e^{-\xi\overline{G}(\mathbf{e}_q)-\alpha
S(\mathbf{e}_q)
};J(\mathbf{e}_q)\right]=\mathbf{E}\left[e^{-\xi\overline{G}(\mathbf{e}_q)-\alpha
S(\mathbf{e}_q)
};J(\overline{G}(\mathbf{e}_q))\right]\mathbf{C}\;,
\end{equation}
for
\begin{eqnarray}\mathbf{C}&=&\lim_{n\to\infty}{\rm diag}\left(\mathbb{P}_i\left(\mathbf{e}_q<\tau_{1/n}^+\right)^{-1}\right)
\mathbf{P}\left(\mathbf{e}_q<\tau^+_{1/n};J(\mathbf{e}_q)\right)\nonumber\\&=&\lim_{n\to\infty}
\mathbf{P}\left(J(\mathbf{e}_q)|S(\mathbf{e}_q)<1/n\right)\;.\label{matrixC}
\end{eqnarray}
Now (\ref{th3ii1b}) will follow straightforward from
(\ref{matrixC}) and (\ref{th3ii1}), since by (\ref{up-crossing}),
$$
\mathbf{P}\left(\mathbf{e}_q<\tau^+_{1/n};J(\mathbf{e}_q)\right)=
\mathbf{I}(q)-\mathbf{P}\left(\mathbf{e}_q>\tau^+_{1/n};J(\mathbf{e}_q)\right)=
\left(\mathbf{I}-e^{-\mathbf{\Xi}(q,0)1/n}\right) \mathbf{I}(q)
$$
and hence
\begin{equation}\label{AandP}\mathbf{C}
=\left[{\rm diag}
\left(\mathbf{\Xi}(q,0)\mathbf{I}\mathbf{e}^{\text
T}\right)\right]^{-1}\mathbf{\Xi}(q,0)\mathbf{I}(q)
.\end{equation} The identity (\ref{th3ii2b}) follows from
(\ref{th3ii1b}) and the Theorem \ref{WH.main}(i). Finally, from
the proof of the Theorem \ref{WH.main}(i) it follows that
\begin{eqnarray} \lefteqn{\mathbf{E}\left[
e^{\alpha I(\mathbf{e}_q)-\xi \underline{G}(\mathbf{e}_q)};
J(\mathbf{e}_q) \right]}\nonumber
\\&&=\mathbf{E}\left[
e^{\alpha I(\mathbf{e}_q)-\xi \underline{G}(\mathbf{e}_q)};
J(\underline{G}(\mathbf{e}_q)) \right]
\mathbf{\Delta}_{\mathbf{\pi}}^{-1}\widehat{\mathbf{P}}\left(
J(\overline{G}(\mathbf{e}_q))\right)^{\text T}
\mathbf{\Delta}_{\mathbf{\pi}}, \nonumber\end{eqnarray} which
completes the proof of (\ref{th3ii2}) in a view of
(\ref{th3ii1b}).

\subsection{Proof of Theorem \ref{rogozin}}\label{rogozinproof}

For a general matrix $\mathbf{A}$ with the distinct eigenvalues $\lambda_i$ (hence with the
independent eigenvectors $\mathbf{s}_i$),
under assumption that
$q\mathbf{I}+\mathbf{A}$ is invertible, using Frullani
integral and the representation $\mathbf{A}=\mathbf{S}\;{\rm diag}\{\lambda_i\}\;\mathbf{S}^{-1}$
with $\mathbf{S}=(\mathbf{s}_1,\ldots, \mathbf{s}_N)$, we can derive the following identity:
\begin{equation}\label{frullani}q\left(q\mathbf{I}+\mathbf{A}\right)^{-1}=
\exp\left\{\int_0^\infty
(e^{-\mathbf{A}x}-\mathbf{I})\frac{1}{x}e^{-q
x}\;dx\right\}.\end{equation}

\begin{lemma}\label{firstdecomposition}
Under assumption (\ref{commute}) for $\xi\geq 0$,
\begin{eqnarray*} \lefteqn{ \mathbf{E}\left[e^{i\alpha
X(\mathbf{e}_q)-\xi
\mathbf{e}_q};J(\mathbf{e}_q)\right]}\\&&=\exp\left\{\int_0^\infty\int_{[0,\infty)}
\left(\exp\left\{-\xi  t +i\alpha x
\right\}-1\right)\frac{1}{t}e^{-qt}\;\mathbf{P}(X(t)\in dx; J(t))
\;dt\right\} \nonumber\\&&\ \ \
\cdot\exp\left\{\int_0^\infty\int_{(-\infty,0)}
\left(\exp\left\{-\xi t +i\alpha x
\right\}-1\right)\frac{1}{t}e^{-qt}\;\mathbf{P}(X(t)\in dx; J(t))
\;dt\right\} \;.
\end{eqnarray*}
\end{lemma}
\proof By additivity of process $X(t)$ there exists a matrix
$\mathbf{F}$ such that $\mathbf{E}\exp\{i\alpha
X(t)\}=\exp\{\mathbf{F}(i\alpha)t\}$ (see Prop. XI.2.2, p. 311 of
Asmussen (2003)). Note that this matrix also has distinct
eigenvalues. Then we have:
\begin{eqnarray*}
\lefteqn{ \mathbf{E}\left[e^{i\alpha X(\mathbf{e}_q)-\xi
\mathbf{e}_q};J(\mathbf{e}_q)\right]=\int_0^\infty qe^{-qt}
\exp\left\{-(\xi \mathbf{I}-\mathbf{F}(i\alpha))t
\right\}\;dt}\nonumber\\&&= \exp\left\{\int_0^\infty
\left(\exp\left\{-(\xi\mathbf{I}-\mathbf{F}(i\alpha))t\right\}-\mathbf{I}\right)\frac{1}{t}e^{-qt}\;dt\right\}
\nonumber\\&&=\exp\left\{\int_0^\infty\int_{\mathbf{R}}
\left(\exp\left\{-\xi  t +i\alpha x
\right\}-1\right)\frac{1}{t}e^{-qt}\;\mathbf{P}(X(t)\in dx; J(t))
\;dt\right\}.
\end{eqnarray*}
Note that by identity (\ref{commute}) the matrices
$$\int_0^\infty\int_{[0,\infty)}
\left(\exp\left\{-\xi  t +i\alpha x
\right\}-1\right)\frac{1}{t}e^{-qt}\;\mathbf{P}(X(t)\in dx; J(t))
\;dt$$ and
$$\int_0^\infty\int_{(-\infty,0)}
\left(\exp\left\{-\xi  t +i\alpha x
\right\}-1\right)\frac{1}{t}e^{-qt}\;\mathbf{P}(X(t)\in dx; J(t))
\;dt$$ commutes.
This gives the assertion of the lemma by the factorization.
\hspace*{\fill}\mbox{$\square$}

From Lemma \ref{firstdecomposition} and Theorem \ref{WH.main}
we have
\begin{eqnarray}
&&\exp\left\{\int_0^\infty\int_{[0,\infty)}
\left(\exp\left\{-\xi  t +i\alpha x
\right\}-1\right)\frac{1}{t}e^{-qt}\;\mathbf{P}(X(t)\in dx; J(t))
\;dt\right\} \nonumber\\&&\ \ \
\cdot\exp\left\{\int_0^\infty\int_{(-\infty,0)}
\left(\exp\left\{-\xi t +i\alpha x
\right\}-1\right)\frac{1}{t}e^{-qt}\;\mathbf{P}(X(t)\in dx; J(t))
\;dt\right\}\nonumber\\&&= \mathbf{H}(\alpha,
\xi)\mathbf{T}(\alpha, \xi), \label{unique}
\end{eqnarray}
where \begin{equation}\label{mathbfH}
\mathbf{H}(\alpha, \xi)= \mathbf{E}\left[ e^{i\alpha
S(\mathbf{e}_q)-\xi \overline{G}(\mathbf{e}_q)};
J(\mathbf{e}_q)\right]\mathbf{I}(q)^{-1}\end{equation}
and \begin{equation}\label{mathbfT}
\mathbf{T}(\alpha, \xi)=\mathbf{\Xi}(q,0)^{-1}{\rm diag}
\left(\mathbf{\Xi}(q,0)\mathbf{e}^{\text
T}\right)\mathbf{\Delta}_{\mathbf{\pi}}^{-1}\widehat{\mathbf{E}}\left[
e^{i\alpha I(\mathbf{e}_q)-\xi \underline{G}(\mathbf{e}_q)};
J(\underline{G}(\mathbf{e}_q))\right]^{\text T}
\mathbf{\Delta}_{\mathbf{\pi}}.
\end{equation}
From Theorem \ref{WH.main}(i) it follows that matrices
$\mathbf{H}(\alpha, \xi)$ and $\mathbf{T}(\alpha, \xi)$ are
invertible. Thus,
\begin{eqnarray}
\lefteqn{
\mathbf{H}^{-1}(\alpha,
\xi)\exp\left\{\int_0^\infty\int_{[0,\infty)}
\left(\exp\left\{-\xi  t +i\alpha x
\right\}-1\right)\right.}\nonumber\\&&\hspace{5cm}\left.\frac{1}{t}e^{-qt}\;\mathbf{P}(X(t)\in dx; J(t))
\;dt\right\} \nonumber\\&&
=\mathbf{T}(\alpha,
\xi)\exp\left\{-\int_0^\infty\int_{(-\infty,0)}
\left(\exp\left\{-\xi t +i\alpha x
\right\}-1\right)\right.\nonumber\\&&\hspace{5cm}\left.\frac{1}{t}e^{-qt}\;\mathbf{P}(X(t)\in dx; J(t))
\;dt\right\}. \label{unique2}
\end{eqnarray}
Moreover, each entry of the matrix $\mathbf{H}(\alpha, \xi)$ is
analytical in the upper half of the complex plane. The same
concerns also the matrix $\mathbf{H}^{-1}(\alpha, \xi)$. Thus each
entry of the LHS of (\ref{unique2}) extends analytically to the
lower half of the complex plane in $\alpha$ and similarly each
entry of the matrix on the RHS of (\ref{unique2}) extends
analytically to the upper half of the complex plane in $\alpha$.
Hence matrices on both sides of (\ref{unique2}) can be defined in
the whole $\alpha$-plane. Observe that each entry of
these matrices is bounded function. Indeed, from definitions (\ref{mathbfH})
and (\ref{mathbfT})
by Jensen inequality it follows that
each entry of matrices $\mathbf{H}(\alpha, \xi)$
and $\mathbf{T}(\alpha, \xi)$ is bounded in respective regions.
Note that reciprocal of
determinant of $\mathbf{H}(\alpha, \xi)$ is also bounded within any bounded circle. Thus on any circle each entry of $\mathbf{H}^{-1}(\alpha, \xi)$ is bounded.
Similarly, one can prove that each entry of the second factors of the RHS and LHS of (\ref{unique2})
is bounded.
Thus by Liouville's Theorem
each entry of (\ref{unique2}) must be a constant. Putting $\alpha=\xi=0$ gives the
assertion of the theorem.

\subsection{Proof of Theorem \ref{kendall}}\label{kendallproof}
We will use now the martingale technique to prove Theorem \ref{kendall}
(see Borovkov and Burq (2001) for similar considerations in the
case of spectrally negative L\'{e}vy processes).

\begin{lemma}\label{burq1}
We have, $$ \int_y^\infty
\frac{dxe^{-\Phi(q)x}}{\kappa^\prime(\Phi(q))}= \int_0^\infty
te^{-qt} \int_y^\infty \frac{dx}{x}\sum_{j=1}^N
\frac{\mathbf{h}_{j}(\Phi(q))}{\mathbf{h}_i(\Phi(q))}
\mathbb{P}_i(\tau_x^+\in dt; J(t)=j).
$$
\end{lemma}
\proof From (\ref{alaGirsanov}) we have
$$e^{-\Phi(q)x}=\mathbb{E}_i
e^{-q\tau_x^+}\mathbf{h}_{J(\tau_x^+)}(\Phi(q))/\mathbf{h}_{i}(\Phi(q))\;.$$
Differentiating with respect to $q$ and noting that
$\Phi^\prime(q)=1/\kappa^\prime(\Phi(q))$ we have
\begin{equation}\label{firststep}
xe^{-\Phi(q)x}\frac{1}{\kappa^\prime(\Phi(q))}=\int_0^\infty
te^{-qt} \sum_{j=1}^N
\frac{\mathbf{h}_{J(t)}(\Phi(q))}{\mathbf{h}_i(\Phi(q))}
\mathbb{P}_i(\tau_x^+\in dt; J(t)=j)\;.
\end{equation}
The proof completes by dividing left-hand side and right-hand side
of (\ref{firststep}) by $x$ and integrating them w.r.t. $dx$ over
$(y,\infty)$. \hspace*{\fill}\mbox{$\square$}
\begin{lemma}\label{burq2}
We have,
$$\int_y^\infty
\frac{dxe^{-\Phi(q)x}}{\kappa^\prime(\Phi(q))}= \int_0^\infty
e^{-qt} \int_y^\infty
\sum_{j=1}^N\frac{\mathbf{h}_j(\Phi(q))}{\mathbf{h}_i(\Phi(q))}
\mathbb{P}_i\left(X(t)\in dx; J(t)=j\right)\;dt.$$
\end{lemma}
\proof From Corollary XI.2.6 of Asmussen (2003) we have
$$a=\kappa^\prime(\Phi(q))\mathbb{E}_i^{\Phi(q)}\tau_a^+
+\mathbf{h}_i(\Phi(q))-\mathbb{E}_i^{\Phi(q)}\mathbf{h}_{J(\tau_x^+)}(\Phi(q))$$ and
hence
\begin{equation}\label{atoinfty}\lim_{a\to\infty}
\frac{\mathbb{E}_i^{\Phi(q)}\tau_a^+}{a}=\frac{1}{\kappa^\prime(\Phi(q))}\;.\end{equation}
Let $T_A=\int_0^\infty \mathbf{I}_A(X(t))\;dt$ be the time spend
by our process in the set $A$. Note that
\begin{equation}\label{twosidedeq}
\tau_a^+-T_{(-\infty,0]}\leq T_{(0,a]}\leq
\tau_a^++T^{J(\tau_a^+)}_{(-\infty,0]}\;,\end{equation} where
$T^{J(\tau_a^+)}_{(-\infty,0]}$ denotes the time spend in
$(-\infty,0]$ by the process
$\{(X(t+\tau^+_a)-a,J(\tau_a^++t)),t\ge 0\}$. Moreover,
\begin{eqnarray*}
\lefteqn{\max_i\mathbb{E}_i^{\Phi(q)}T_{(-\infty,0]} =\max_i
\mathbb{E}_i^{\Phi(q)}\int_0^\infty\mathbf{I}_{(-\infty,0]}(X(t))\;dt}\\&&\leq
\max_i\mathbb{E}_i^{\Phi(q)}\int_0^\infty\mathbf{I}_{(-\infty,0]}(X(t))\;dt\\&&\leq
\max_i\mathbb{E}_i^{\Phi(q)}\int_0^\infty e^{-qX(t)}\;dt\leq
\max_i\int_0^\infty\mathbb{E}_i^{\Phi(q)} e^{-qX(t)}\;dt\\&&\leq
\max_{i,j}\frac{\mathbf{h}_{j}(\Phi(q))}{\mathbf{h}_i(\Phi(q))}
\int_0^\infty e^{-qt}\;dt=
\frac{1}{q}\max_{i,j}\frac{\mathbf{h}_{j}(\Phi(q))}{\mathbf{h}_i(\Phi(q))}<\infty
\;.
\end{eqnarray*}
Taking expectation from both sides of (\ref{twosidedeq}) and
dividing by $a$ we derive:
\begin{equation}\label{secondstep}
\lim_{a\to\infty}
\frac{\mathbb{E}_i^{\Phi(q)}T_{(0,a]}}{a}=\frac{1}{\kappa^\prime(\Phi(q))}\;.\end{equation}
Observe that for $0<a<b<c<\infty$
$$T_{(a,b]}+T_{(b,c]}=T_{(a,c]}$$
and hence
$$\mathbb{E}_i^{\Phi(q)}T_{(a,b]}=c(b-a)$$
which together with (\ref{secondstep}) gives:
\begin{equation}\label{firdstep}
\mathbb{E}_i^{\Phi(q)}T_{(0,a]}=\frac{a}{\kappa^\prime(\Phi(q))}\;.
\end{equation}
Above can be rewritten in the following way:
\begin{eqnarray*}
\lefteqn{\frac{a}{\kappa^\prime(\Phi(q))}=\mathbb{E}_i^{\Phi(q)}T_{(0,a]}=
\mathbb{E}_i^{\Phi(q)}\int_0^\infty\mathbf{I}_{(0,a]}(X(t))\;dt}\\&&=
\int_0^\infty\mathbb{P}_i^{\Phi(q)}\left(X(t)\in
(0,a]\right)\;dt\\&&= \int_0^\infty
dt\int_{(0,a]}e^{\Phi(q)x-qt}\sum_{j=1}^N\frac{\mathbf{h}_j(\Phi(q))}{\mathbf{h}_i(\Phi(q))}
\mathbb{P}_i\left(X(t)\in dx; J(t)=j\right)\\&&=
\int_{(0,a]}e^{\Phi(q)x} \int_0^\infty e^{-qt}
\sum_{j=1}^N\frac{\mathbf{h}_j(\Phi(q))}{\mathbf{h}_i(\Phi(q))}
\mathbb{P}_i\left(X(t)\in dx; J(t)=j\right)\;dt\;.
\end{eqnarray*}
Since this equation holds for any $a>0$, we have
$$
\frac{dx}{\kappa^\prime(\Phi(q))}=e^{\Phi(q)x} \int_0^\infty
e^{-qt} \sum_{j=1}^N\frac{\mathbf{h}_j(\Phi(q))}{\mathbf{h}_i(\Phi(q))}
\mathbb{P}_i\left(X(t)\in dx; J(t)=j\right)\;dt
$$
which completes the proof by integrating them w.r.t. $dx$ over
$(y,\infty)$. \hspace*{\fill}\mbox{$\square$}

From Lemma \ref{burq1} and \ref{burq2} we have the following
equality of Laplace transforms: \begin{eqnarray*} \lefteqn{
\int_0^\infty te^{-qt} \int_y^\infty \frac{dx}{x}\sum_{j=1}^N
\frac{\mathbf{h}_{j}(\Phi(q))}{\mathbf{h}_i(\Phi(q))}
\mathbb{P}_i(\tau_x^+\in dt; J(t)=j)}\\&&= \int_0^\infty e^{-qt}
\int_y^\infty
\sum_{j=1}^N\frac{\mathbf{h}_j(\Phi(q))}{\mathbf{h}_i(\Phi(q))}
\mathbb{P}_i\left(X(t)\in dx; J(t)=j\right)\;dt
\end{eqnarray*}
that implies the equality of the following measures:
\begin{equation}\label{newstep}
\sum_{j=1}^N\mathbf{h}_{j}(\Phi(q))\left(t\frac{dx}{x}\mathbb{P}_i(\tau_x^+\in
dt; J(t)=j)-\mathbb{P}_i\left(X(t)\in dx; J(t)=j\right)dt\right)=0
\end{equation}
for each $i=1,\ldots,N$. This is equivalent to:
$$\left(t dx\mathbf{P}(\tau_x^+\in
dt;J(t))-x\mathbf{P}\left(X(t)\in dx;
J(t)\right)dt\right)\mathbf{h}(\Phi(q))=\mathbf{0}\;.
$$
Choosing now the independent vectors $\mathbf{h}(\Phi(q_1)), \mathbf{h}(\Phi(q_2)),\ldots, \mathbf{h}(\Phi(q_N))$ completes the proof.

\subsection{Proof of Theorem \ref{ballot}}\label{ballotproof}
By Kendall's identity given in Theorem \ref{kendall} it suffices
to prove that
$$\mathbf{P}(X(t)\in dx,I(t)=0;J(t))dt=\frac{1}{c}\mathbf{P}(\tau_x^+\in dt; J(t))dx$$
or that for all $q>0$ and sufficiently large $s>0$:
\begin{eqnarray}
\lefteqn{q\int_0^\infty e^{-q t}dt\int_0^\infty e^{s
x}\mathbf{P}(X(t)\in dx,I(t)=0;J(t))}\nonumber\\&&=
\frac{q}{c}\int_0^\infty e^{-q t}\int_0^\infty e^{sx}dx
\mathbf{P}(\tau_x^+\in dt; J(t))\label{startingpoint1}
\end{eqnarray}
that is equivalent to
\begin{equation}\label{startingpoint2}
\lim_{\alpha \to \infty}\mathbf{E}\left[e^{sX(\mathbf{e}_q)+\alpha
I(\mathbf{e}_q)};J(\mathbf{e}_q)\right]= \frac{q}{c} \int_0^\infty
e^{sx} \mathbf{E}\left[e^{-q\tau_x^+}; J(\tau_x^+)\right]dx\;.
\end{equation}
We prove (\ref{startingpoint2}) passing from its left-hand side to
its right-hand side. Let $\overline{q}=q-\kappa(s)$. The change of
measure (\ref{alaGirsanov}) and Wiener-Hopf factorization given in
Theorem \ref{WH.main}(ii) yields:
\begin{eqnarray*}
\lefteqn{\lim_{\alpha \to
\infty}\mathbf{E}\left[e^{sX(\mathbf{e}_q)+\alpha
I(\mathbf{e}_q)};J(\mathbf{e}_q)\right]}\\&&=\lim_{\alpha \to
\infty}\mathbf{\Delta}_{\mathbf{h}}(s)\mathbf{E}^{s}\left[e^{\alpha
I(\mathbf{e}_q)-\kappa(s)\mathbf{e}_q};J(\mathbf{e}_q)\right]\mathbf{\Delta}_{\mathbf{h}}(s)^{-1}\\&&
=\lim_{\alpha \to
\infty}\frac{q}{\overline{q}}\mathbf{\Delta}_{\mathbf{h}}(s)\mathbf{E}^{s}\left[e^{\alpha
I(\mathbf{e}_{\overline{q}})};J(\mathbf{e}_{\overline{q}})\right]\mathbf{\Delta}_{\mathbf{h}}(s)^{-1}\\&&
=\lim_{\alpha \to \infty}\mathbf{\Delta}_{\mathbf{h}}(s)
q\left(\frac{\overline{q}}{\alpha}\mathbf{I}-\mathbf{F}_s(\alpha)/\alpha\right)^{-1}
\mathbf{\Delta}_{\mathbf{\pi}_s}^{-1}\frac{1}{\alpha}\widehat{\mathbf{\Xi}}_s(\overline{q},-\alpha)^{\text
T}\left[\widehat{\mathbf{\Xi}}_s(\overline{q},0)^{-1}\right]^{\text
T}\mathbf{\Delta}_{\mathbf{\pi}_s}\mathbf{\Delta}_{\mathbf{h}}(s)^{-1}.
\end{eqnarray*}
Note that
$$\lim_{\alpha \to
\infty}\left(\frac{\overline{q}}{\alpha}\mathbf{I}-\mathbf{F}_s(\alpha)/\alpha\right)^{-1}=
\left(-c\mathbf{I}\right)^{-1}=-\frac{1}{c}\mathbf{I}
$$
and from (\ref{Xi})
$$\lim_{\alpha \to
\infty}\frac{1}{\alpha}\widehat{\mathbf{\Xi}}_s(\overline{q},-\alpha)^{\text
T}=-\mathbf{I}.$$ Hence $$\lim_{\alpha \to
\infty}\mathbf{E}\left[e^{sX(\mathbf{e}_{q})+\alpha
I(\mathbf{e}_{q})};J(\mathbf{e}_{q})\right]= \frac{q}{c}
\mathbf{\Delta}_{\mathbf{h}}(s)\mathbf{\Delta}_{\mathbf{\pi}_s}^{-1}
\left[\widehat{\mathbf{\Xi}}_s(\overline{q},0)^{-1}\right]^{\text
T}
\mathbf{\Delta}_{\mathbf{\pi}_s}\mathbf{\Delta}_{\mathbf{h}}(s)^{-1}
\;.
$$
Using classical arguments for reversed process we can proceed as
follows:
\begin{eqnarray*}
\lefteqn{\frac{q}{c}
\mathbf{\Delta}_{\mathbf{h}}(s)\mathbf{\Delta}_{\mathbf{\pi}_s}^{-1}
\left[\widehat{\mathbf{\Xi}}_s(\overline{q},0)^{-1}\right]^{\text
T} \mathbf{\Delta}_{\mathbf{\pi}_s}
\mathbf{\Delta}_{\mathbf{h}}(s)^{-1}}\\&&= \frac{q}{c}
\mathbf{\Delta}_{\mathbf{h}}(s)\int_0^\infty
\mathbf{\Delta}_{\mathbf{\pi}_s}^{-1}
\widehat{\mathbf{E}}^s\left[e^{-\overline{q}\tau_x^+};J(\tau_x^+)\right]^{\text
T}\mathbf{\Delta}_{\mathbf{\pi}_s}\;dx
\mathbf{\Delta}_{\mathbf{h}}(s)^{-1}
\\&&= \frac{q}{c}\int_0^\infty
\mathbf{\Delta}_{\mathbf{h}}(s)\mathbf{E}^s\left[e^{-\overline{q}\tau_x^+};J(\tau_x^+)\right]
\mathbf{\Delta}_{\mathbf{h}}(s)^{-1}\;dx
\\&&=
\frac{q}{c}\int_0^\infty
\mathbf{E}\left[e^{-(\overline{q}+\kappa(s))\tau_x^+
+sx};J(\tau_x^+)\right] \;dx\;,
\end{eqnarray*}
which is the right-hand side of (\ref{startingpoint2}).

\section*{Acknowledgements}
ZP would like to thank Andreas Kyprianou who 
works with him at the beginning of this project and who gave lots of
valuable comments.   
This work is partially supported by the Ministry of Science
and Higher Education of Poland under the grant N N2014079 33 (2007-2009).

\end{document}